\documentclass [10pt,a4paper,draft] {article}
\usepackage [applemac] {inputenc}
\usepackage[french] {babel}
\usepackage {amsfonts}
\usepackage {amsmath}
\usepackage {amssymb}

\begin{document}

\begin{center}
\begin{LARGE}
Notes sur l'indice des alg\`ebres de Lie (II)
\end{LARGE}
\bigskip

par : Mustapha RA\"IS 
\footnotemark \footnotetext[1]{Universit\'e de POITIERS -
D\'epartement de Math\'ematiques - T\'el\'eport 2, Boulevard Marie et Pierre Curie - BP
30179 -
86962 FUTUROSCOPE CHASSENEUIL Cedex} 
\end{center}

\bigskip

\vspace{15mm}

Ce texte est une suite \`a : ``\textsl{Notes sur l'indice des alg\`ebres de Lie (I)}'',
dont les notations sont
conserv\'ees pour l'essentiel.

%%%%%%%%%%%%%%%%%%%%%%%
\section{Les champs de vecteurs invariants et leurs d\'eriv\'ees}

\noindent
$\centerdot$ Soient $\mathfrak{g}$ une alg\`ebre de Lie (de dimension finie sur un
corps $k$, disons
$k = \mathbb{R}$ ou $\mathbb{C}$) et $P : \mathfrak{g} \longrightarrow \mathfrak{g}$
une
application polynomiale (i.e. un champ de vecteurs sur $\mathfrak{g}$, polynomial),
homog\`ene de
degr\'e $m \geq 1$, et invariante sous l'action du groupe adjoint $G$ de
$\mathfrak{g}$. On a donc :
$$
                (1)\qquad         P(\hbox{Ad}(g)x) = \hbox{Ad}(g)P(x)\quad (x \in \mathfrak{g},\ g \in G).
$$
La forme infinit\'esimale de cette propri\'et\'e est :
$$
        dP(x).[y,x] = [y,P(x)]\quad (x\ \hbox{et}\ y\ \hbox{dans}\ \mathfrak{g}).
$$
(Ici et plus loin, on utilise les notations habituelles du calcul diff\'erentiel ; par
exemple : $dP(x)$ est la
valeur au point $x$ de la d\'eriv\'ee premi\`ere de la fonction $P$, c'est donc un
endomorphisme de
l'espace vectoriel $\mathfrak{g}$, et dans le premier membre de l'\'egalit\'e ci-dessus,
on applique cet
endomorphisme au ``vecteur'' $[y,x]$.)

\bigskip
        On notera que, de cette \'egalit\'e, il r\'esulte imm\'ediatement que $P(x)$ appartient au
centre $z(z(x))$
du centralisateur $z(x)$ de $x$.

\vskip 7mm
\noindent
$\centerdot$ On \'ecrit la formule de Taylor pour $P$ :
$$
        (2)\qquad                 P(x+ty) = \sum_{0\leq k \leq m}\ {t^k\over k!}\, d^kP(x).y^{(k)}\quad (x\ 
\hbox{et}\ y \ \hbox{dans}\ \mathfrak{g}, \ \ t \ \hbox{dans} \ k).
$$
Chaque terme $d^kP(x).y^{(k)}$, consid\'er\'e comme une fonction de $x$ et de $y$, est
polynomial
homog\`ene de degr\'e $(m-k)$ en $x$, de degr\'e $k$ en $y$ et on a (une formule d'\'echange) :
$$
        (3)\qquad                {1\over k!}\ d^kP(x).y^{(k)} = {1\over (m-k)!}\ d^{m-k}P(y).x^{(m-k)}.
$$
En particulier :
$$
        m! \,  P(x) = \mathcal{P}.x^{(m)}\quad (x\in \mathfrak{g})
$$
o\`u $\mathcal{P} = d^mP$ est la d\'eriv\'ee \`a l'ordre ``maximum'' $m$ (c'est donc une
$m$-forme
sym\'etrique sur $\mathfrak{g}$, de degr\'e z\'ero en $x$, donc ind\'ependante de $x$).

\vskip 7mm
\noindent
$\centerdot$ L'invariance de $P$ se propage en l'invariance des diverses fonctions
intervenant dans
le second membre de la formule (2). On a en fait (pour tous $x,y$ dans
$\mathfrak{g}$, et $0 \leq k
\leq m$) : 
$$
        (4)\qquad                [z, d^kP(x).y^{(k)}] = d^{k+1}P(x).[z,x].y^{(k)} + k\,
d^kP(x).[z,y].y^{(k-1)}.
$$

\vskip 7mm
\noindent
$\centerdot$ Soit $(h,e,f)$ un $sl(2)$-triplet inclus dans $\mathfrak{g}$ (s'il en
existe). En utilisant
les formules (4), on trouve (seule $[h,e] = 2e$ intervient) :
$$
        (5)\qquad                [h,\ d^k\, P(h).e^{(k)}] = 2k\ d^k\, P(h).e^{(k)}
$$
$$
        (6)\qquad                [e,\ d^k\, P(h).e^{(k)}] = -2\ d^{k+1}\, P(h).e^{(k+1)}.
$$

\bigskip
        Comme : $(ade)^{m-k}d^kP(h).e^{(k)} = (-2)^{m-k}d^mP(h).e^{(m)} = (-2)^{m-k}\, m!\,
P(e)$, on
voit que sous l'hypoth\`ese : $P(e) \not = 0$, les vecteurs $d^kP(h).e^{(k)}$\,
$(0\leq k \leq m)$ sont
lin\'eairement ind\'ependants (ce sont des vecteurs propres de $adh$, associ\'es
respectivement
aux valeurs propres $0,2,4,\ldots ,2m$).

\vskip 12mm

%%%%%%%%%%%%%%%%%%%%%%%
\section{Le cas d'une alg\`ebre de Lie simple}

        Dor\'enavant, $\mathfrak{g}$ est une alg\`ebre de Lie simple et le corps de base est le
corps des
complexes. Soient $p_1, p_2,\ldots ,p_r$ un syst\`eme de g\'en\'erateurs homog\`enes, de
degr\'es respectifs
$m_1+1,\ldots , m_r+1$, alg\'ebriquement ind\'ependants, de l'alg\`ebre
$\mathbb{C}[\mathfrak{g}]^G$
des fonctions polyn\^omes $G$-invariantes sur $\mathfrak{g}$. Ainsi $r$ est le rang de
$\mathfrak{g}$
et $m_1, m_2,\ldots , m_r$ sont les exposants de $\mathfrak{g}$. Pour chaque
entier\break$j=
1,2,\ldots ,r$ on note : $P_j : \mathfrak{g} \longrightarrow \mathfrak{g}$ le
gradient de $p_j$, calcul\'e
au moyen de la forme de Killing $B$ de $\mathfrak{g}$ :
$$
        ({d\over dt})_0\ p_j(x+ty) = \ <dp_j(x),y>\ = B(P_j(x),y)\quad (x \ \hbox{et}\ y\
\hbox{dans}\
\mathfrak{g}).
$$
Chaque $P_j$ est un champ de vecteurs invariant, polynomial homog\`ene de degr\'e $m_j$,
et on peut
appliquer le paragraphe 1 \`a chaque $P_j$, et \`a un \'el\'ement nilpotent
\underline{r\'egulier} $e$, de
sorte que $(h,e,f)$ est un $sl(2)$-triplet principal. On a donc : $[h,P_j(e)] =
2m_j\, P_j(e)$, et
$[e,P_j(e)] = 0$. Les vecteurs $P_j(e)\ (1\leq j \leq r)$ sont lin\'eairement
ind\'ependants (d'apr\`es un
ancien r\'esultat de Kostant) et sont des vecteurs primitifs pour la repr\'esentation de
$\mathfrak{a} =
\mathbb{C}h + \mathbb{C}e + \mathbb{C}f$ dans $\mathfrak{g}$ (restriction \`a
$\mathfrak{a}$ de
la repr\'esentation adjointe de $\mathfrak{g}$).

\vskip 7mm
\noindent
$\centerdot$ \textbf{Notes} : Posons $v_{j,k} = d^kP_j(h).e^{(k)}\quad (1\leq j \leq
r,\ \ 0\leq k \leq
m_j)$. On a donc : 
$$
        (7)\qquad        [h, v_{j,k}] = 2k\, v_{j,k}
$$
$$
        (8)\qquad        [e, v_{j,k}] = -2\, v_{j,k+1}
$$
avec la convention : $v_{j,m_j+1} = 0$.

\bigskip
        Posons $w_{j,k} = d^kP_j(h).f^{(k)}\ (1\leq j \leq r$ et $0 \leq k \leq m_j)$. Par
les m\^emes calculs
que plus haut ($[h,f] = -2f$ intervenant \`a la place de $[h,e] = 2e$), on trouve : 
$$
        [h, w_{j,k}] = -2k\, w_{j,k}\quad \hbox{et : }\ [f,w_{j,k}] = 2w_{j,k+1}.
$$
Il se trouve que l'ensemble des vecteurs $(v_{j,k})_{1\leq j \leq r,\ 0\leq k \leq
m_j}$ conjointement
avec les $w_{j,k}$ $(1\leq j \leq r, \ 1 \leq k \leq m_j)$ forme une base de
l'espace vectoriel
$\mathfrak{g}$, et qu'on a la d\'ecomposition triangulaire : $\mathfrak{g} =
\mathfrak{n}_- \oplus
\mathfrak{h} \oplus \mathfrak{n}_+$, avec $\displaystyle \mathfrak{h} = \sum_{1\leq
j \leq r}\
\mathbb{C}\, P_j(h),$\quad  et 
$\mathfrak{n}_+ = \displaystyle \sum_{\begin{subarray}{1}
1\leq j \leq r\\
1\leq k \leq m_j
\end{subarray}} \mathbb{C}\ v_{j,k}$ et 
$\mathfrak{n}_- = \displaystyle \sum_{\begin{subarray}{1}
1\leq j \leq r\\
1\leq k \leq m_j
\end{subarray}} \mathbb{C}\ w_{j,k}$. Clairement, la $h$-graduation de
$\mathfrak{g}$ est mise en
\'evidence : $\mathfrak{g}^{(0)} = \mathfrak{h},\ \mathfrak{n}_+ = \displaystyle
\sum_{\ell > 0}\,
\mathfrak{g}^{(\ell)}$ et $\mathfrak{n}_- = \displaystyle \sum_{\ell < 0}\,
\mathfrak{g}^{(\ell)}$. Le lecteur int\'eress\'e pourra trouver la d\'emonstration de ces
faits dans \cite{Ra}.

\vskip 7mm
\noindent
$\centerdot$ Consid\'erons l'ensemble des fonctions $\varphi_t$ sur $\mathfrak{g}$
d\'efinies par : pour
tout $x \in \mathfrak{g},\ \varphi_t(x) = \varphi(x + th)$ o\`u $\varphi$ d\'ecrit
l'ensemble
$\mathbb{C}[\mathfrak{g}]^\mathfrak{g}$ des fonctions polyn\^omes invariantes et $t$
d\'ecrit
$\mathbb{C}$. Il est bien connu et imm\'ediat que deux telles fonctions sont en
involution relativement
au crochet de Lie-Poisson sur $\mathfrak{g}$. Ce qui pr\'ec\`ede montre que l'espace
engendr\'e par les
$[e, \nabla \varphi_t(e)]$ est de dimension \'egale \`a la moiti\'e de la dimension de
l'orbite nilpotente de
$e$. Ceci s'interpr\`ete en terme de syst\`emes compl\`etement int\'egrables (voir par
exemple \cite{Bol}).

\vskip 12mm

%%%%%%%%%%%%%%%%%%%%%%%
\section{Le normalisateur du centralisateur d'un \'el\'ement nilpotent}

\bigskip
On reprend les notations du paragraphe pr\'ec\'edent, \`a ceci pr\`es que l'\'el\'ement
nilpotent $e$ n'est plus
forc\'ement un \'el\'ement r\'egulier. On fera toutefois l'hypoth\`ese suivante :

\medskip
        ``Le centre $\delta(e)$ du centralisateur $z(e)$ de $e$ est engendr\'e par la famille
$(P_j(e))_{1\leq
j \leq r}$''.

\vskip 5mm
Soient $j_1,\ldots , j_s$ les entiers tels que $(P_{j_1}(e),\ldots , P_{j_s}(e))$
soit         une base de
$\delta(e)$. Pour simplifier les notations, on posera :
$$
        Q_1 = P_{j_1},\ldots , Q_s = P_{j_s},\ m'_1 = m_{j_1},\ldots ,m'_s = m_{j_s}, \
\delta(e) = \delta\
\hbox{et}\ z(e) = z
$$
et on supposera que les entiers $m'_1,\ldots ,m'_s$ sont rang\'es dans l'ordre
croissant (avec la
terminologie de \cite{Ri}, $(m'_1,\ldots ,m'_s)$ est la suite des exposants de
$(\mathfrak{g}, e))$.

\vskip 5mm
        On pose : $y_j = dQ_j(e).h\quad (1\leq j \leq s)$.

\vskip 7mm
\noindent
\textbf{3.1. Lemme :}~\textit{Soit $\eta$ le normalisateur de $z$ dans
$\mathfrak{g}$. On a :}
$$
        \eta = z \oplus \sum_{1\leq j \leq s}\, \mathbb{C}y_j.
$$

\vskip 7mm
\noindent
\textbf{D\'emonstration} :  Compte-tenu des formules (4), on a :
$$
        [e, y_j] = -2m'_j\, Q_j(e).
$$
Ceci prouve que les $y_j\ (1\leq j \leq s)$ sont lin\'eairement ind\'ependants, et qu'en
notant $V$
l'espace vectoriel engendr\'e par les $y_j \ (1\leq j \leq s)$, $ad(e)$ induit une
bijection de $V$ sur
$\delta$. On sait par ailleurs (voir par exemple \cite{Tau} 17.5.12) que $\dim \eta
= \dim z + \dim
\delta$. D'o\`u le r\'esultat.

\vskip 7mm
\noindent
\textsc{Remarque} : Posons $z_j = Q_j(e)\ (1\leq j \leq s)$. Toujours avec l'aide
des formules (4), on
voit que :
$$
        [f, z_j] = dQ_j(e).[f,e] = -y_j.
$$
        Ceci est une autre d\'emonstration du lemme, compte-tenu de \cite{Tau} (17.5.6).

\vskip 7mm
%\vfill\eject
 %%%%
\noindent
\textbf{3.2.}~On a, compte-tenu des formules (7) :
$$
        [h, y_j] = 2(m'_j - 1)y_j \quad (1\leq j \leq s)
$$
$$
        [h, z_j] = 2m'_j z_j \quad (1\leq j \leq s).
$$

\bigskip
Dans la suite, on s'int\'eresse au calcul des $[y_i, z_j]\quad (1\leq i,j \leq s)$.
Imm\'ediatement,
compte-tenu du fait que $[y_i, z_j]$ est de $h$-graduation $2(m'_i + m'_j -1)$, on
d\'eduit : 
$$
        [y_i, z_j] = 0
$$
lorsque $(m'_i + m'_j -1)$ n'est pas un \underline{exposant} de $(\mathfrak{g},e)$,
et en particulier
lorsque $m'_i + m'_j > m'_s + 1$.

\bigskip
        On va retrouver ce r\'esultat en \'etablissant un lien avec la notion de convolution ou
d\'eplacement des
invariants (\cite{A}, \cite{G}).

\medskip
        On pose : $\omega_{ij} = B(Q_i, Q_j)$ o\`u, comme indiqu\'e plus haut, $B$ est la forme
de Killing de
$\mathfrak{g}$ (on peut d'ailleurs remplacer $B$ par n'importe quel multiple
scalaire non nul de $B$).

\vskip 7mm
 %%%%
\noindent
\textbf{3.3. Lemme : }

\medskip
$        (1)\quad        [y_i, z_j] = [y_j, z_i]$

\medskip
$        (2)\quad        [y_i, z_j] = 2m'_j\ dQ_i(e).Q_j(e)$

\medskip
$        (3)\quad        dQ_i(e).Q_j(e) = (L_j\, Q_i)(e)$

\smallskip
\qquad \textit{o\`u $L_j$ est l'op\'erateur de d\'erivation le long du champ de vecteurs
$Q_j$.}

\medskip
$        (4)\quad        [y_i, z_j] = {m'_im'_j\over m'_i+m'_j}\ \nabla \omega_{ij}(e)$

\smallskip
\qquad \textit{o\`u $\nabla \omega_{ij}$ est le gradient de la fonction $\omega_{ij}$.}

\vskip 7mm
\noindent
\textbf{D\'emonstration} : (1) Ceci est bien connu (\cite{Pa}) :
$$
        0 = [f,[z_i,z_j]], \ \hbox{et}\ [f,z_k] = -y_k
$$

\vskip 5mm
\noindent
(2) \ $[z_j,y_i] = [z_j, dQ_i(e).h] = d^2Q_i(e).[z_j,e].h + dQ_i(e).[z_j,h]$.

\bigskip
        Comme $[z_j, e] = 0$ (puisque $z_j \in \delta$) et $[z_j,h] = -2m'_j z_j$, on
arrive \`a :
$$
        [z_j,y_i] = -2m'_j\ dQ_i(e).Q_j(e).
$$

\vskip 5mm
\noindent
(3) \ Par d\'efinition de l'op\'erateur diff\'erentiel $L_j$ : 

\begin{eqnarray*}
        L_j Q_i(x) &=        &({d\over dt})_0\ Q_i(x + t\, Q_j(x))\\
                                                                                                                                        &=         &dQ_i(x).Q_j(x)
\end{eqnarray*}

\vskip 5mm
\noindent
(4) \         Avec $x$ et $y$ dans $\mathfrak{g}$, on a :
$$
        <d\omega_{ij}(x), y>\ = B(dQ_i(x).y,\ Q_j(x)) + B(Q_i(x),\, dQ_j(x).y)
$$
        \begin{eqnarray*}
        B(dQ_i(x).y,\ Q_j(x))         &=                &d^2 q_i(x) (y,Q_j(x))\\
                                                                                                                                                                                                                                                                                        &=                &d^2q_i(x) (Q_j(x),y)\\
                                                                                                                                                                                                                                                                                        &=                &B(dQ_i(x).Q_j(x),y)
        \end{eqnarray*}

o\`u $q_i = p_{j_i}$. Donc :
$$
        <d\omega_{ij}(x),y>\ = B(dQ_i(x).Q_j(x) + dQ_j(x).Q_i(x),y)
$$
et
$$
        \nabla \omega_{ij}(x) = dQ_i(x).Q_j(x) + dQ_j(x).Q_i(x)
$$
D'o\`u : 
        \begin{eqnarray*}
        m'_i\, m'_j\ \nabla \omega_{ij}(e) &= &m'_i (m'_j\, dQ_i(e).Q_j(e)) + m'_j (m'_i\,
dQ_j(e).Q_i(e))\\
        &= &(m'_i + m'_j)({1\over 2}\, [y_i, z_j] + {1\over 2}\, [y_j, z_i])\\
 &= &(m'_i + m'_j) [y_i, z_j]
\end{eqnarray*}

\medskip
et enfin :  
$$
        [y_i, z_j] = {m'_i m'_j\over m'_i + m'_j}\ \nabla \omega_{ij}(e)
$$

\vskip 7mm
 %%%%
\noindent
\textbf{3.4.}~On voit appara\^\i tre les ``produits scalaires'' $\omega_{ij}$ d\'ej\`a
pr\'esents ailleurs, en
particulier dans les travaux d'Arnold et Givental (\cite{A}, \cite{G}) sous le nom
de convolution ou
d\'eplacement des invariants.

\vskip 7mm
\noindent
$\centerdot$ Chaque $\omega_{ij}$ est une fonction polyn\^ome invariante sur
$\mathfrak{g}$,
homog\`ene de degr\'e : $(m'_i + m'_j)$.  Notons $I_+$ l'id\'eal de
$\mathbb{C}[\mathfrak{g}]^\mathfrak{g}$ engendr\'e par $p_1, p_2,\ldots , p_r$. Il
existe des
constantes $c^k_{ij}\ (1\leq k \leq r)$, bien d\'etermin\'ees, telles que :
$$
        \omega_{ij} = \sum_{1\leq k \leq r}\ c^k_{ij}\, p_k\quad \mod I^2_+
$$
La fonction $\displaystyle \omega'_{ij} = \sum_{1\leq k \leq r}\ c^k_{ij}\, p_k$
s'appelle la ``partie
lin\'eaire'' de $\omega_{ij}$ dans les travaux d'Arnold-Givental d\'ej\`a cit\'es. On a alors :
$$
        \nabla \omega_{ij}(e) = \nabla \omega'_{ij}(e) = \sum_{1 \leq k \leq r}\ c^k_{ij}\,
P_k(e)
$$
car le gradient d'une fonction appartenant \`a $I^2_+$ est nul en tout nilpotent de
$\mathfrak{g}$.

\vskip 7mm
 %%%%
\noindent
\textbf{3.5.}~Supposons dor\'enavant, en plus de l'hypoth\`ese d\'ej\`a faite, que les
exposants de
$\mathfrak{g}$ soient 2 \`a 2 distincts : $1 = m_1 < m_2 <\cdots < m_r$. Alors, comme
indiqu\'e dans
\cite{Ri}, une base de $\delta$ est constitu\'ee par les $P_j(e)$ \underline{qui sont
non nuls}.
Autrement dit, avec les notations pr\'ec\'edentes, $\{j_1, j_2,\ldots ,j_s \}$ est
l'ensemble des indices $j$
tels que $P_j(e) \not = 0$. Dans cette circonstance, $\nabla \omega_{ij}(e)$ est non
nul si et seulement
si : $\displaystyle \omega'_{ij} = \sum_{1\leq k \leq s}\, \alpha^k_{ij}\, q_k$, o\`u
les $\alpha^k_{ij}$
ne sont pas tous nuls.

\vskip 5mm
\noindent
$\centerdot$ \textbf{Conclusion} : 
        \begin{enumerate} 
        \item Lorsque $m'_i + m'_j -1$ n'est pas un exposant de $(\mathfrak{g},e),\
\omega'_{ij} = 0$ et :
$[y_i,z_j] = 0$. C'est le cas en particulier lorsque $m'_i + m'_j -1 > m'_s$.

        \medskip
        \item Lorsque $m'_i + m'_j = 1+ m'_k$, pour un entier $k$ alors bien d\'etermin\'e, on
a : $\omega'_{ij}
= \alpha_i q_k$, $\alpha_i$ \'etant un nombre complexe, et $[y_i, z_j] = \beta_i
Q_k(e)$, o\`u $\beta_i$
est un nombre complexe qui est un multiple de $\alpha_i$ : 
$$
        \beta_i = {m'_i m'_j\over 1+m'_k}\, \alpha_i
$$
de sorte que $[y_i,z_j] \not = 0$ ssi $\omega'_{ij} \not = 0$.
        \end{enumerate}

\vskip 7mm
 %%%%
\noindent
\textbf{3.6.}~Soit $\mathcal{A}$ la matrice $([y_i, z_j])_{1\leq i,j \leq s}$. Cette
matrice
est pseudo-triangulaire :
$$
        [y_i, z_j] = 0 \quad \hbox{lorsque}\ i+j > s+1
$$
et : $[y_i, z_{s+1-i}] = \beta_i\, Q_s(e)$. La matrice $\mathcal{A}$ \'etant
consid\'er\'ee comme une
matrice \`a coefficients dans l'alg\`ebre sym\'etrique de $\mathfrak{g}$, on calcule son
d\'eterminant : 
$$
        \det \mathcal{A} = \beta_1\, \beta_2\cdots \beta_s(Q_s(e))^s.
$$

\noindent
$\centerdot$ Soit $\mathcal{B} = (\omega'_{ij})_{1\leq i,j\leq s}$. C'est une matrice
pseudo-triangulaire \`a coefficients fonctions polyn\^omes sur $\mathfrak{g}$ et 
$$
        \det \mathcal{B} = \alpha_1 \alpha_2 \cdots \alpha_s(q_s)^s.
$$

\vskip 5mm
\indent
Il vient la \underline{conclusion} : $\hbox{ind}(\eta, \delta) = 0$ ssi $\det
\mathcal{A} \not = 0$ ssi
$\det \mathcal{B} \not = 0$. Plus pr\'ecis\'ement, on rappelle que : $\hbox{ind}(\eta,
\delta) = \dim
\delta - \, \hbox{rg}(\mathcal{A})$, o\`u $\hbox{rg}(\mathcal{A})$ est le rang de $A$.

\vskip 7mm
 %%%%
\noindent
\textbf{3.7.}~Revenons au cas particulier o\`u $e$ est un \'el\'ement nilpotent
\underline{r\'egulier}. Dans ce
cas, on a : $\hbox{ind}(\eta, z) = 0$ d'apr\`es un th\'eor\`eme de Panyushev (\cite{Pa},
5.6), i.e. :
$$
        \det \mathcal{A} = \gamma(P_r(e))^r, \quad \gamma \in \mathbb{C}^*
$$
et $P_r(e)$ est l'\'el\'ement de plus grande $h$-graduation dans $z : [h, P_r(e)] =
2m_r\, P_r(e)$.

\bigskip
        Par ailleurs, le fait que $\det \mathcal{B}$ soit non nul appara\^\i t d\'ej\`a dans
\cite{A} et principalement
dans \cite{G}. Dans ce dernier article cit\'e, on trouvera les calculs explicites pour
les alg\`ebres de type
$B_n, C_n, D_n, F_4$ et $E_6$, des produits scalaires $\omega_{ij}$ (pour un choix
particulier des
$p_j$). On peut utiliser ces calculs pour d\'eterminer les nombres $\hbox{ind}(\eta,
\delta)$ \`a condition
que l'hypoth\`ese faite : ``$\displaystyle \delta(e) = \sum_{1\leq j \leq r}\,
\mathbb{C}\, P_j(e)$'' soit
v\'erifi\'ee. Pour les alg\`ebres de Lie classiques, c'est le cas pour tous les \'el\'ements
nilpotents des alg\`ebres
$sl(n),\, so(2n+1), \, sp(2n)$ et pour certains types de nilpotents de $so(2n)$ ;
dans toutes ces
situations, il a \'et\'e prouv\'e par Panyushev (\cite{Pa}, theorem 4.7) que le groupe $N$
associ\'e \`a $\eta$
admet un nombre fini d'orbites dans $\eta^*$, et en particulier que :
$\hbox{ind}(\eta, \delta) = 0$. Le
point de vue adopt\'e ici (passage par les $\nabla \omega_{ij}(e)$) n'apporte rien de
nouveau. Toutefois,
les calculs explicites de Givental (\cite{G}) des $\omega_{ij}$ pour des alg\`ebres
exceptionnelles, par
exemple $F_4$ et $E_6$, permettent d'\'ecrire la matrice $\mathcal{A} = ([y_i,
z_j])_{i,j}$ pour les
nilpotents particuliers v\'erifiant l'hypoth\`ese rappel\'ee plus haut, et par suite de
calculer
$\hbox{ind}(\eta, \delta)$.

\bigskip
        Il reste au moins \`a faire la liste des \'el\'ements nilpotents auxquels on peut
appliquer cette m\'ethode
 et \`a expliquer l'intervention des $\omega_{ij}$.

\vskip 12mm

\renewcommand{\refname}{Bibliographie} %%%pour m\'emoire : \cite{ } dans le texte

\end {document}